\input amstex.tex
\documentstyle{amsppt}
\magnification=\magstep1
%
%
\def\section#1{\par\bigpagebreak\par\bigpagebreak
    \csname subhead\endcsname #1\endsubhead\par\bigpagebreak}
\def\subsection#1{\par\medpagebreak
    \csname subsubhead\endcsname #1 \endsubsubhead}
\def\Theorem#1#2{\csname proclaim\endcsname{Theorem #1} #2 
    \endproclaim}
\def\Corollary#1{\csname proclaim\endcsname{Corollary} #1 \endproclaim}
\def\Proposition#1#2{\csname proclaim\endcsname{Proposition #1} #2 
    \endproclaim}
\def\Lemma#1#2{\csname proclaim\endcsname{Lemma #1} #2 \endproclaim}
\def\Remark#1#2{\remark{Remark {\rm #1}} #2 \endremark}
\def\Proof#1{\demo{Proof} #1\qed\enddemo}
\def\Proofof#1#2{\demo{Proof of #1} #2\qed\enddemo}
%
%
\define\ep{\epsilon}
\define\dt{\delta}
\define\Dt{\Delta}
\define\ld{\lambda}

\define\Z{{\Bbb Z}}
\define\Q{{\Bbb Q}}

\define\K{{\Bbb K}}

\define\scpr#1#2{\langle{#1},\,{#2}\rangle}

\define\bb#1{\text{\bf #1}}

\leftheadtext\nofrills{A.N.\,Kirillov and M.\,Noumi}
\rightheadtext\nofrills{$q$-Difference raising operators}
\topmatter
\title\nofrills
$q$-Difference raising operators for Macdonald polynomials \\
and the integrality of transition coefficients
\endtitle
\author
Anatol N.\, KIRILLOV
\footnote"$^{\ast1}$"
{\hbox{Department of Mathematical Sciences, University of Tokyo {\it and}}
\newline
\hbox{\quad St.\,Petersburg Branch of the Steklov Mathematical Institute}}
{\rm and} 
Masatoshi NOUMI
\footnote"$^{\ast2}$"
{\hbox{Department of Mathematics, Kobe University}}
\endauthor
\address
A.N.\,Kirillov: 
\newline\indent
{\eightpoint
Department of Mathematical Sciences, University of Tokyo,
\newline\indent
Komaba, Meguro, Tokyo 153, Japan;
\newline\indent
Steklov Mathematical Institute, 
\newline\indent
Fontanka 27, St.\,Petersburg 191011, Russia
}
\newline\indent
M.\,Noumi: 
\newline\indent
{\eightpoint
Department of Mathematics, Faculty of Science, Kobe University,
\newline\indent
Rokko, Kobe 657, Japan}
\endaddress
\email 
kirillov\@ker.c.u-tokyo.ac.jp\quad
noumi\@math.s.kobe-u.ac.jp
\endemail
\endtopmatter
\document
\section{Introduction}
The purpose of this paper is study certain $q$-difference raising operators 
for Macdonald polynomials (of type $A_{n-1}$) which are originated from the 
$q$-difference - reflection operators introduced in our previous 
paper \cite{KN}. 
These operators can be regarded as a $q$-difference version 
of the raising operators for Jack polynomials introduced by 
L.\,Lapointe and L.\, Vinet \cite{LV1, LV2}. 
As an application of our $q$-difference raising operators, we will give a 
proof of the integrality of the double Kostka coefficients which had been 
conjectured by I.G.\,Macdonald \cite{Ma}, Chapter VI.
We will also determine their {\it quasi-classical limits}, 
which give rise to (differential) raising operators for Jack polynomials. 
(See also {\it Notes} at the end of Introduction.) 
\par\medpagebreak
Let $\K=\Q(q,t)$ be the field of rational functions in two indeterminates 
$(q,t)$ 
and $\K[x]^W$ the algebra of symmetric polynomials in 
$n$ variables $x=(x_1,\cdots,x_n)$ over $\K$, $W$ being the symmetric group 
$\frak{S}_n$ of degree $n$. 
The {\it Macdonald polynomials} $P_\ld(x;q,t)$ 
are a family of symmetric polynomials parametrized by partitions, 
and they form a $\K$-basis of $\K[x]^W$.  
They are characterized as the joint eigenfunctions in $\K[x]^W$ for 
Macdonald's commuting family of $q$-difference operators
$$
D_r=
t^\binom{r}{2}\sum\Sb I\subset[1,n] \\ |I|=r\endSb
\prod\Sb i\in I\\ j\not\in I\endSb \frac{tx_i-x_j}{x_i-x_j} 
\prod_{i\in I}T_{q,x_i}
\quad (r=0,1,\cdots, n).  
\tag{1}
$$
For each partition $\ld$, 
the Macdonald polynomial $P_\ld(x;q,t)$ is the unique joint 
eigenfunction of $D_r$ ($r=0,1,\cdots,n$) that has the leading term 
$m_\ld(x)$ under the dominance order of partitions when it is expressed 
as a linear combination of monomial symmetric functions $m_\mu(x)$.  
As in \cite{Ma} (VI.8.3), we also use another normalization 
$J_\ld(x;q,t)=c_\ld(q,t) P_\ld(x;q,t)$,
called the ``integral form'' of $P_\ld(x;q,t)$. 
\par
We now define the two kinds of $q$-difference operators $K^{+}_m$ 
and $K^{-}_m$ ($m=0,1,\ldots,n$) as follows:
$$
\align
K^+_m
&=\sum\Sb J\subset[1,n]\\|J|=m\endSb 
\prod_{j\in J} x_j
\sum\Sb I\subset J\endSb
(-t^{m-n+1})^{|I|} t^{\binom{|I|}{2}}
\prod\Sb i\in I\\j\in[1,n]\backslash I\endSb 
\frac{tx_i-x_j}{x_i-x_j}
\prod_{i\in I} T_{q,x_i},
\tag{2}\\
K^-_m&=
\sum\Sb J\subset[1,n]\\|J|=m\endSb 
\prod_{j\in J} x_j
\sum\Sb I\subset J\endSb
(-t)^{m-|I|} t^{\binom{m-|I|}{2}}
\prod\Sb i\in I\\j\in[1,n]\backslash I\endSb 
\frac{x_i-tx_j}{x_i-x_j}
\prod_{j\in[1,n]\backslash I} T_{q,x_j}. 
\endalign
$$
These operators are $W$-invariant and preserve the ring of symmetric 
polynomials
$\K[x]^W$.
\Theorem{A}{
For each $m=0,1,\cdots,n$, 
the $q$-difference operator $K_m=K^+_m$ (resp. $K^-_m$) is a raising operator 
for Macdonald polynomials $J_\ld(x;q,t)$ in the sense that 
$$
K_m J_\ld(x;q,t)=J_{\ld+(1^m)}(x;q,t)
\tag{3}
$$
for any partition $\ld$ with $\ell(\ld)\le m$. 
}
\noindent
(See Theorem 2.2 in Section 2.) 
\par
Theorem A implies that, 
for any partition $\ld=(\ld_1,\cdots,\ld_n)$, 
the Macdonald polynomial $J_\ld(x;q,t)$ 
is obtained by a successive application of the operators $K_m$ starting 
from $J_{0}(x;q,t)=1$:
$$
J_{\ld}(x;q,t)=(K_n)^{\ld_n}(K_{n-1})^{\ld_{n-1}-\ld_n}\cdots
(K_1)^{\ld_1-\ld_2}(1).  
\tag{4}
$$
{}From this expression, we can show that $J_{\ld}(x;q,t)$ is a linear 
combination 
of monomial symmetric functions with coefficients in $\Z[q,t]$.  
Furthermore we have 
\Theorem{B}{
For any partition $\ld$ and $\mu$,  
the double Kostka coefficient $K_{\ld,\mu}(q,t)$ is a polynomial in 
$q$ and $t$ with integral coefficients. 
} 
\noindent
(See Theorem 2.4.)  
Theorem B gives a partial affirmative answer 
to the conjecture of Macdonald proposed in 
\cite{Ma}, (VI.8.18?) (apart from the positivity of the coefficients). 
\par\medpagebreak
After some preliminaries on Macdonald's $q$-difference operators 
$D_r$,  we formulate our main results in Section 2 and show how 
Theorem A implies Theorem B. 
We will propose in Section 3 some determinantal formulas related 
to our raising operators $K^\pm_m$.  
The proof of Theorem A will be given in Section 4 by analyzing the action 
of the operators $K^\pm_m$ on the generating function of Macdonald polynomials. 
In Section 5, we will include a similar construction of lowering operators
for Macdonald polynomials. 
\par\medpagebreak
{\eightpoint 
{\it Notes\,}: 
In \cite{KN}, we constructed the raising operators for Macdonald polynomials 
by means of the Dunkl operators due to I.\,Cherednik, in an analogous way 
as L.\,Lapointe and L.\,Vinet \cite{LV1, LV2} did for Jack polynomials.  
We also gave an application to the integrality of transition coefficients 
of Macdonald polynomials. 
Although these raising operators involve reflection operators, 
they act as $q$-difference operators on symmetric functions;
the explicit forms of the corresponding $q$-difference operators 
are also determined in \cite{KN}. 
After this work, we found a direct, elementary proof of the fact that
the $q$-difference operators in question have the property of raising 
operators for Macdonald polynomials. 
In this paper we begin with introducing those $q$-difference operators 
directly, and show in an elementary way (without affine Hecke algebras 
and Dunkl operators) that they are the raising operators for Macdonald 
polynomials that we want. 
In view of its elementary nature, we decided to make the present 
paper as self-contained as possible, and also to repeat here 
the proof of integrality of transition coefficients for the sake of 
reference. 
For this reason, this paper has some intersection with our previous paper 
\cite{KN} (Section 1 and some part of Section 2).  
{}From the viewpoint of this paper, our previous paper could be thought of 
as explaining the meaning of the $q$-difference raising operators 
in relation to affine Hecke algebras and Dunkl operators.  
}
\par\newpage
\section{\S1:  Macdonald's $q$-difference operators}
%
In this section, we will make a brief review of some basic properties 
of the Macdonald polynomials (associated with the root system of type 
$A_{n-1}$, or the symmetric functions with two parameters) and 
the commuting family of $q$-difference operators
which have Macdonald polynomials as joint eigenfunctions.  
For details, see Macdonald's book \cite{Ma}. 
\par\smallpagebreak\par
Let $\K=\Q(q,t)$ be the field of rational functions in two indeterminates 
$q$, $t$ 
and consider the ring $\K[x]=\K[x_1,\cdots,x_n]$ of polynomials in $n$ 
variables $x=(x_1,\cdots,x_n)$ with coefficients in $\K$.  
Under the natural action of the symmetric group $W=\frak{S}_n$ of 
degree $n$, the subring of all symmetric polynomials will be denoted 
by $\K[x]^{W}$. 
\par 
The {\it Macdonald polynomials} $P_{\ld}(x)=P_{\ld}(x;q,t)$ 
(associated with the root system of type $A_{n-1}$) 
are symmetric polynomials parametrized by the 
{\it partitions } $\ld=(\ld_1,\cdots,\ld_n)$
($\ld_i\in\Z$, $\ld_1\ge\cdots\ge\ld_n\ge0$).  
They form a $\K$-basis of the invariant ring $\K[x]^W$
and are characterized as the joint eigenfunctions 
of a commuting family of $q$-difference operators 
$\{D_r\}_{r=0}^n$. 
For each $r=0,1,\cdots,n$, the $q$-difference operator $D_r$ 
is defined by 
$$
D_r=t^{\binom{r}{2}} \sum\Sb I\subset[1,n] \\ |I|=r \endSb 
\prod\Sb i\in I\\ j\not\in I\endSb 
\frac{tx_i-x_j}{x_i-x_j} \, \prod_{i\in I}T_{q,x_i}, 
\tag{1.1}
$$
where 
$T_{q,x_i}$ stands for the $q$-shift operator in the variable 
$x_i$ : 
$(T_{q,x_i}f)(x_1,\cdots,x_n)=f(x_1,\cdots,q x_i,\cdots,x_n)$. 
The summation in (1.1) is taken over all subsets $I$ of the interval 
$[1,n]=\{1,2,\cdots,n\}$ consisting of $r$ elements. 
Note that $D_0=1$ and 
$D_n=t^{\binom{n}{2}}T_{q,x_1}\cdots T_{q,x_n}$. 
Introducing a parameter $u$, we will use the generating function
$$
D_x(u)=\sum_{r=0}^n (-u)^r D_r
\tag{1.2} 
$$
of these operators $\{D_r\}_{r=0}^n$. 
Note that the operator $D_x(u)$ has the determinantal expression
$$
D_x(u)=\frac{1}{\Dt(x)}\det(x_j^{n-i}(1-ut^{n-i}T_{q,x_j}))_{1\le i,j\le n},
\tag{1.3}
$$
where $\Dt(x)=\prod_{1\le i<j\le n}(x_i-x_j)$ is the difference product 
of $x_1,\cdots,x_n$. 
It is well known that the $q$-difference operators $D_r$ ($0\le r\le n$) 
commute with each other, or equivalently, $[D_x(u), D_x(v)]=0$. 
Furthermore the Macdonald polynomial $P_{\ld}(x)$ satisfies the 
$q$-difference equation
$$
D_x(u) P_{\ld}(x)=d^n_\ld(u) P_{\ld}(x), \quad\text{with}\quad 
d^n_\ld(u)=\prod_{i=1}^n (1-ut^{n-i}q^{\ld_i}), 
\tag{1.4}
$$
for each partition $\ld=(\ld_1,\cdots,\ld_n)$. 
Recall that each $P_\ld(x)$ can be written in the 
form
$$
P_\ld(x)=m_{\ld}(x) + \sum_{\mu<\ld} u_{\ld\mu} m_{\mu}(x)\quad(u_{\ld\mu}
\in \K),
\tag{1.5}
$$
where, for each partition $\mu$, $m_{\mu}(x)$ stands for the monomial 
symmetric function of type $\mu$, and $\le$ is the dominance 
order of partitions. 
The Macdonald polynomials $P_\ld(x)$ are determined uniquely by the conditions 
(1.4) and (1.5). 
The ``integral form'' $J_\ld(x)=J_\ld(x;q,t)$ is defined by the normalization 
$$
J_\ld(x)=c_\ld P_\ld(x) ,\quad c_\ld=\prod_{s\in\ld} 
(1-t^{\ell(s)+1}q^{a(s)}), 
\tag{1.6}
$$
where for each square $s=(i,j)$ in the diagram of a partition $\lambda$,
the numbers $\ell(s)=\ld'_j-i$ and $a(s)=\ld_i-j$ are respectively the
{\it leg-length}  and the {\it arm-length} of $\ld$ at $s\in\ld$
(\cite{Ma}, (VI.6.19)).
\par
We recall now the ``reproducing kernel'' for the Macdonald polynomials. 
Consider another set of variables $y=(y_1,\cdots,y_m)$ and assume that 
$m\le n$.  
We define the function $\Pi(x,y)=\Pi(x,y;q,t)$ by
$$
\Pi(x,y)=\prod\Sb i\in[1,n]\\ j\in[1,m]\endSb
\frac{(tx_iy_j;q)_\infty}{(x_iy_j;q)_\infty}, 
\tag{1.7}
$$
where $(x;q)_\infty=\prod_{k=0}^\infty (1-xq^k)$. 
The convergence of the infinite product above may be 
understood in the sense of formal power series.
It is known that the function $\Pi(x,y)$ has the expression 
$$
\Pi(x,y)=\sum_{\ell(\ld)\le m} b_\ld P_\ld(x) P_\ld(y) \quad(b_\ld\in \K),
\tag{1.8}
$$
where the summation is taken over all partitions $\ld$ with length $\le m$, 
and each partition $\ld=(\ld_1,\cdots,\ld_m,0,\cdots,0)$ with 
$\ell(\ld)\le m$
is identified with the truncation $(\ld_1,\cdots,\ld_m)$ when it is used as 
the suffix for $P_\ld(y)$. 
The coefficients $b_\ld$ in (1.8) are determined as
$$
b_\ld=\prod_{s\in\ld} 
\frac{1-t^{\ell(s)+1}q^{a(s)}}{1-t^{\ell(s)}q^{a(s)+1}}.
\tag{1.9}
$$
\par
We remark that, by (1.4), expression (1.8) is equivalent to the formula
$$
D_x(u)\Pi(x,y)=(u;t)_{n-m} D_y(u t^{n-m}) \Pi(x,y).
\tag{1.10}
$$
Since 
$$
D_x(u)=\sum_{I\subset[1,n]} (-u)^{|I|} t^{\binom{|I|}{2}}
\prod\Sb i\in I\\ j\not\in I\endSb
\frac{tx_i-x_j}{x_i-x_j} \, \prod_{i\in I}T_{q,x_i}, 
\tag{1.11}
$$
We have 
$$
D_x(u) \Pi(x,y) = \Pi(x,y) F(u;x,y),
\tag{1.12}
$$
where
$$
F(u;x,y)=\sum_{I\subset[1,n]} (-u)^{|I|} t^{\binom{|I|}{2}}
\prod\Sb i\in I\\ j\not\in I\endSb 
\frac{tx_i-x_j}{x_i-x_j} \, 
\prod\Sb i\in I\\ k\in[1,m]\endSb
\frac{1-x_iy_k}{1-tx_iy_k}.
\tag{1.13}
$$
Hence formula (1.10) is also equivalent to 
$$
F(u;x,y)=(u;t)_{n-m} F(ut^{n-m};y,x). 
\tag{1.14}
$$
(See also \cite{MN}.)
%
\section{\S2: $q$-Difference raising operators and transition coefficients} 
We now define the $q$-difference operators $K^+_m$ and $K^-_m$ 
($m=0,1,\ldots,n$) as follows:  
$$
\align
K^+_m
&=\sum\Sb J\subset[1,n]\\|J|=m\endSb 
\prod_{j\in J} x_j
\sum\Sb I\subset J\endSb
(-t^{m-n+1})^{|I|} t^{\binom{|I|}{2}}
\prod\Sb i\in I\\j\in[1,n]\backslash I\endSb 
\frac{tx_i-x_j}{x_i-x_j}
\prod_{i\in I} T_{q,x_i},
\tag{2.1}\\
K^-_m&=
\sum\Sb J\subset[1,n]\\|J|=m\endSb 
\prod_{j\in J} x_j
\sum\Sb I\subset J\endSb
(-t)^{m-|I|} t^{\binom{m-|I|}{2}}
\prod\Sb i\in I\\j\in[1,n]\backslash I\endSb 
\frac{x_i-tx_j}{x_i-x_j}
\prod_{j\in[1,n]\backslash I} T_{q,x_j}. 
\endalign
$$
These operators are {\it dual} to each other in the sense that 
$$
K^-_m=(-t)^m t^{\binom{m}{2}} (K^+_m)^\vee \, T_{q,x_1}\cdots T_{q,x_n}. 
\tag{2.2}
$$
where the superscript $(\cdot)^\vee$ stands for the involution on 
$q$-difference 
operators induced from the transformation $q\to q^{-1}$, $t\to t^{-1}$ 
(cf. \cite{Ma}, (VI.8.5)).
If $m=0$ or $m=n$, the operators $K^\pm_m$ reduce to
$$
K^+_0=1,\  K^+_n=x_1\cdots x_nD_x(t)\quad\text{and}\quad
K^-_0=T_{q,x_1}\cdots T_{q,x_n},\  K^-_n=x_1\cdots x_nD_{x}(t). 
\tag{2.3}
$$
\par
We remark first that these $q$-difference operators are $W$-invariant 
since their definition do not depend on the ordering of the coordinates 
$x_1,\ldots,x_n$.  
Hence they transform symmetric functions to symmetric functions. 
Note also that the coefficients of the operator $\Dt(x) K^+_m$
are in $\Z[t^{\pm 1}][x]$ and that those of $\Dt(x) K^-_m$ are in $\Z[t][x]$. 
If $R$ is any subring of $\K=\Q(q,t)$ and $f(x)$ is an alternating polynomial 
in $R[x]$, then the symmetric polynomial $\Dt(x)^{-1}f(x)$ also 
have coefficients in $R$, i.e. $\Dt(x)^{-1}f(x)\in R[x]^W$.  
{}From these remarks, we obtain 
\Proposition{2.1}{
The $q$-difference operators $K^{\pm}_m$ $(m=0,1,\cdots,n)$
preserve the ring of symmetric polynomials $\K[x]^W$.  
Furthermore, the operators $K^-_m$ preserve the ring $\Z[q,t][x]^W$ 
of symmetric polynomials with coefficients in $\Z[q,t]$. 
}
We now state our main theorem.
\Theorem{2.2}{
For any partition $\ld$ with $\ell(\ld)\le m$ $(m=0,\ldots,n)$, we have 
$$
K^+_mJ_\ld(x)=J_{\ld+(1^m)}(x)\quad\text{and}
\quad K^-_mJ_\ld(x)=J_{\ld+(1^m)}(x).
\tag{2.4}
$$
}
\noindent
Note that $\ld+(1^m)=(\ld_1+1,\ldots,\ld_m+1,0,\ldots,0)$ 
if $\ell(\ld)\le m$. 
When $\ld=0$, we have $J_{(1^m)}(x)=e_m(x) (t;t)_m$, 
$e_m(x)$ being the elementary symmetric function of degree $m$.
Hence, formulas in (2.4) imply the following: 
$$
\align
&\sum\Sb J\subset[1,n]\\|J|=m\endSb 
\prod_{j\in J} x_j
\sum\Sb I\subset J\endSb
(-t^{m-n+1})^{|I|} t^{\binom{|I|}{2}}
\prod\Sb i\in I\\j\notin I\endSb 
\frac{tx_i-x_j}{x_i-x_j}
=e_m(x) (t;t)_m,
\tag{2.5}\\
&\sum\Sb J\subset[1,n]\\|J|=m\endSb 
\prod_{j\in J} x_j
\sum\Sb I\subset J\endSb
(-t)^{m-|I|} t^{\binom{m-|I|}{2}}
\prod\Sb i\in I\\j\notin I\endSb 
\frac{x_i-tx_j}{x_i-x_j}
=e_m(x) (t;t)_m.
\endalign
$$
Theorem 2.2 will be proved later in Section 4.
In this section, we will discuss some consequences of our theorem.  
\par\medpagebreak
{}From Theorem 2.2, it follows that the Macdonald polynomial $J_\ld(x)$ 
for any partition $\ld$ can be obtained from the constant function $J_0(x)=1$ 
by a successive application of the raising operators $K_m=K^+_m$ or $K^-_m$. 
Namely we have
$$
J_\ld(x)
=(K_n)^{\ld_n}(K_{n-1})^{\ld_{n-1}-\ld_n}\cdots (K_1)^{\ld_1-\ld_2}(1), 
\tag{2.6}
$$
or equivalently
$$
J_\ld(x)
=K_{\mu_1} K_{\mu_2}\cdots K_{\mu_s}(1), 
\tag{2.7}
$$
where $\mu=(\mu_1,\cdots,\mu_s)$ is the conjugate partition $\ld'$ of $\ld$. 
Since the operators $K^-_m$ preserve the ring $\Z[q,t][x]^W$, we have
$J_{\ld}(x)\in\Z[q,t]^W$.  
In other words,
\Theorem{2.3}{
For any partition $\ld$, the Macdonald polynomial 
$J_\ld(x)$ is expressed as a linear combination of monomial symmetric 
functions $m_\mu(x)$ with coefficients in $\Z[q,t]$. 
}
\noindent
Note that this statement is valid for any $n$.  
Hence it implies that the transition coefficients between $J_{\ld}(x)$ and 
the monomial symmetric functions {\it in infinite variables} are 
polynomials in $q$ and $t$ with integral coefficients. 
\par
{}From this theorem, it is not difficult to conclude the integrality of 
the so-called double Kostka coefficients. 
Recall that the {\it double Kostka coefficients} (or $(q,t)$-Kostka 
coefficients)
$K_{\ld,\mu}(q,t)$ are defined as the transition coefficients between 
the integral forms of Macdonald polynomials and the {\it big Schur functions}: 
$$
J_{\mu}(x;q,t)=\sum_{\ld} K_{\ld,\mu}(q,t) S_{\ld}(x;t),
\tag{2.8}
$$
while Theorem 2.3 is concerned with the transition coefficients 
with the monomial symmetric functions. 
To be more precise, (2.8) should be understood as an equality 
in infinite variables. 
For the definition of $S_\ld(x;t)$, we refer to\cite{Ma}, (III.4.5). 
\Theorem{2.4}{
For any partitions $\ld$ and $\mu$, the double Kostka coefficient
$K_{\ld,\mu}(q,t)$ is a polynomial in $q$ and $t$ with integral coefficients. 
}
\noindent
Theorem 2.4 gives a partial affirmative answer to the conjecture of Macdonald 
proposed in \cite{Ma}, (VI.8.18?), apart from the positivity of the 
coefficients. 
\Proofof{Theorem 2.4}{
As in \cite{Ma}, p.241, the transition coefficients between monomial 
symmetric functions and big Schur functions 
have the form $p(t)/q(t)$, where $p(t), q(t)\in\Z[t]$ and $q(0)=1$. 
Note in particular that they belong to the ring $\Q[t]_{(t)}$ of 
rational functions in $t$, regular at $t=0$. 
Combining this fact with Theorem 2.3, we see that each 
double Kostka coefficients can be written as a finite sum of the form
$$
K_{\ld,\mu}(q,t)=\sum_{k\ge 0} \frac{p^{(k)}_{\ld\mu}(t)}{q^{(k)}_{\ld\mu}(t)}
q^k,
\tag{2.9}
$$
where all $p^{(k)}_{\ld\mu}(t)$ and $q^{(k)}_{\ld\mu}(t)$ belong to $\Z[t]$ 
and $q^{(k)}_{\ld\mu}(0)=1$. 
In particular we have $K_{\ld,\mu}(q,t)\in \Q[t]_{(t)}[q]$.  
We can now apply the duality with respect to $q$ and $t$ again, 
between $J_{\ld}(x;q,t)$ and $J_{\ld'}(x;t,q)$
and between $S_{\ld}(x;t)$ and $S_{\ld'}(x;q)$, to conclude 
$K_{\ld,\mu}(q,t)=K_{\ld',\mu'}(t,q)$ (\cite{Ma},(VI.8.15)).  
Hence, we have 
$K_{\ld,\mu}(q,t)\in \Q[t]_{(t)}[q]\cap  \Q[q]_{(q)}[t]$. 
By using Taylor expansions at $t=q=0$, one can easily see 
that the intersection of the two subalgebras $\Q[t]_{(t)}[q]$ and 
$\Q[q]_{(q)}[t]$ coincides precisely with $\Q[q,t]$.  
Hence we have $K_{\ld,\mu}(q,t)\in \Q[q,t]$. 
In expression (2.9), it means that, for any $k$, 
$p^{(k)}_{\ld\mu}(t)/q^{(k)}_{\ld\mu}(t)\,\in\Q[t]$, 
i.e., $q^{(k)}_{\ld\mu}(t)$ divides $p^{(k)}_{\ld\mu}(t)$. 
Since $q^{(k)}_{\ld\mu}(0)=1$, it follows that 
$p^{(k)}_{\ld\mu}(t)/q^{(k)}_{\ld\mu}(t)\,\in\Z[t]$ 
for all $k$. 
Namely we have $K_{\ld,\mu}(q,t)\in \Z[q,t]$. 
}
\par\medpagebreak
Our raising operators have different expressions, 
from which one can take the {\it quasi-classical} limits as 
$q$ tends to $1$.
\Proposition{2.5}{
For each $m=0,1,\ldots,n$, we have
$$
K^+_m=\frac{1}{\Dt(x)}\sum_{w\in W}\ep(w)w\left(
x_1^{\dt_1}\cdots x_n^{\dt_n}e_{m}
(X_1,\ldots,X_n)\right) 
\tag{2.10}
$$
where $\dt_i=n-i$ $(i=1,\ldots,n)$ and 
$e_m(X_1,\ldots,X_n)$ stands for the $m$-th
elementary symmetric function of the $q$-difference operators 
$$
X_i=x_i(1-t^{m-i+1}T_{q,x_i}) \quad (i=1,\ldots,n).
\tag{2.11}
$$
Similarly we have
$$
K^-_m
=\frac{t^{\binom{n}{2}-\binom{n-m}{2}}}{\Dt(x)}\sum_{w\in W}\ep(w)w\left(
x_1^{\dt_1}\cdots x_n^{\dt_n}e_{m}
(X'_1,\ldots,X'_n)\right) T_{q,x_1}\cdots T_{q,x_n}
\tag{2.12}
$$
where $X'_i= X_it^{-\dt_i}T_{q,x_i}^{-1}$ $(i=1,\cdots,n)$.
}
\noindent
The formulas of Proposition 2.5 will be proved in Section 3 
by using certain {\it determinantal formulas} for the raising operators 
$K^{\pm}_m$.
\par
Let $q=t^{\alpha}$ and let $t\to1$. 
Then the limits of the two $q$-difference operators 
$K^+_m$ and $K^-_m$ give rise to a same 
differential operator:
$$
\align
{\Cal K}_m&=\lim\Sb q=t^\alpha\\ t\to 1\endSb \frac{K^+_m}{(1-t)^m}
=\lim\Sb q=t^\alpha\\ t\to 1\endSb \frac{K^-_m}{(1-t)^m}
\tag{2.13}\\
&=\frac{1}{\Dt(x)}\sum_{w\in W}\ep(w)w\left(
x_1^{\dt_1}\cdots x_n^{\dt_n}e_{m}
({\Cal X}_1,\ldots,{\Cal X}_n)\right),
\endalign
$$
where 
$$
{\Cal X_i}=x_i (\alpha x_i\frac{\partial\ \ }{\partial x_i}+(m-i+1)) 
\quad(i=1,\ldots,n). 
\tag{2.14}
$$
For each partition $\ld$, let $J^{(\alpha)}_\ld(x)$ be the 
Jack polynomial normalized so that the coefficient of 
$m_\ld(x)$ becomes 
$c_\ld(\alpha)=\prod_{s\in\ld}(\alpha \, a(s)+\ell(s)+1)$, 
as in \cite{Ma}, (VI.10.22).
Since 
$$
J^{(\alpha)}_\ld(x)=\lim \Sb t\to 1\endSb 
\frac{J_\ld(x;t^\alpha,t)}{(1-t)^{|\ld|}}, 
\tag{2.15}
$$
Theorem 2.2, combined with (2.13), (2.15), implies the following. 
\Theorem{2.6}{
The differential operators ${\Cal K}_m$ $(m=0,1,\ldots,n)$ in (2.13) are 
raising operators for the Jack polynomials $J^{(\alpha)}_\ld(x)$ 
in the sense that 
$$
{\Cal K}_m J^{(\alpha)}_\ld(x) =J^{(\alpha)}_{\ld+(1^m)}(x)
\tag{2.16}
$$
for any partition $\ld$ with $\ell(\ld)\le m$.
}
\noindent
By a similar argument as we proved Theorem 2.3, we obtain 
\Corollary{ For each partition $\ld$, 
the Jack polynomial $J^{(\alpha)}_\ld(x)$ is expressed 
as a linear combination of monomial symmetric functions 
with coefficients in $\Z[\alpha]$. 
}
Theorem 2.6 can be regarded as a differential version 
(without reflection operators) of the raising operators of 
Lapointe-Vinet \cite{LV1, LV2}. 
%
\section{\S3: Determinantal formulas}
In this section we propose determinantal formulas for some $q$-difference 
operators related to our raising operators. 
The proof of Proposition 2.5 is an application of these determinant formulas. 
Some part of this section will be used also in the proof of Theorem 2.2 in 
Section 4.
\par\medpagebreak
Recall that Macdonald's $q$-difference operator $D_x(u)$ 
has the determinantal formula
$$
D_x(u)=\frac{1}{\Dt(x)}\det(x_j^{n-i}(1-ut^{n-i}T_{q,x_j}))_{1\le i,j\le n}. 
\tag{3.1}
$$
We begin with a comment on this noncommutative determinant. 
Suppose an $n\times n$ matrix $A=(a_{ij})_{i,j}$ with entries 
in a noncommutative algebra has the property
that $a_{ij}a_{k\ell}=a_{k\ell}a_{ij}$ for $i\ne k$ and $j\ne\ell$. 
If this is the case, one has
$$
a_{\sigma(1)1}\cdots a_{\sigma(n)n} 
=a_{1\tau(1)}\cdots a_{n\tau(n)}\quad(\tau=\sigma^{-1})
\tag{3.2}
$$
for any permutation $\sigma$,
since the $n$ elements $a_{\sigma(1)1},\ldots,a_{\sigma(n)n}$
commute with each other. 
Hence we do  not need to distinguish the {\it column} determinant
and the {\it row} determinant:
$$
\det(A)
=\sum_{\sigma\in\frak{S}_n}a_{\sigma(1)1} \ldots a_{\sigma(n)n} 
=\sum_{\sigma\in\frak{S}_n} 
a_{1\sigma(1)} \ldots a_{n\sigma(n)}.
\tag{3.3} 
$$
Since the matrix of (3.1) has this property, and 
there is no ambiguity in the notation. 
We generalize this type of determinantal formula from the viewpoint 
of raising operators. 
\par
In view of the expression (2.1) of our raising operators, let us consider 
the following $q$-difference operators including a parameter $u$: 
$$
\align
K_m(u)
&=\sum\Sb J\subset[1,n]\\|J|=m\endSb x_J
\sum_{I\subset J}(-u)^{|I|} \frac{T_{t,x}^I(\Dt(x))}{\Dt(x)} T_{q,x}^I,
\tag{3.4}\\
L_m(u)
&=\sum\Sb J\subset[1,n]\\|J|=m\endSb x_J
\sum_{I\subset J}(-u)^{m-|I|} \frac{T_{t,x}^{I^{\text{C}}}(\Dt(x))}{\Dt(x)} 
T_{q,x}^{I^{\text{C}}},
\endalign
$$
for each $m=0,1,\cdots,n$, 
where we have used the abbreviation
$$
x_J=\prod_{j\in J},\quad T_{q,x}^I=\prod_{i\in I} T_{q,x_i},
\quad T_{q,x}^{I^{\text{C}}}=\prod_{j\in[1,n]\backslash I} T_{q,x_j}.
\tag{3.5}
$$
Since 
$$
\frac{T_{t,x}^I(\Dt(x))}{\Dt(x)}=t^{\binom{|I|}{2}}
\prod\Sb i\in I\\ j\in[1,n]\backslash I\endSb \frac{tx_i-x_j}{x_i-x_j}, 
\tag{3.6}
$$
the $q$-difference operators $K^\pm_m$ are recovered as the special cases
$$
K^+_m=K_m(t^{m-n+1}),\quad K^-_m=t^{-\binom{n-m}{2}} L_m(t^{m-n+1}). 
\tag{3.7}
$$ 
Introducing another parameter $v$, we set 
$$
K(u,v)=\sum_{m=0}^n v^m K_m(u),\quad L(u,v)=\sum_{m=0}^n v^m L_m(u). 
\tag{3.8}
$$
Then we have the following determinantal formula for $K(u,v)$ and $L(u,v)$. 
\Proposition{3.1}{
$$
\align
K(u,v)&=\frac{1}{\Dt(x)}
\det(x_j^{\dt_i}(1+vx_j(1-ut^{\dt_i}T_{q,x_j})))_{1\le i,j\le n},
\tag{3.9}\\
L(u,v)&=\frac{1}{\Dt(x)}
\det(x_j^{\dt_i}(t^{\dt_i}T_{q,x_j}+vx_j(1-ut^{\dt_i}T_{q,x_j})))_{1\le i,
j\le n},
\endalign
$$
where $\dt_i=n-i$ for $i=1,\ldots,n$. 
When $q=t$, these formulas reduce to
$$
\align
K(u,v)&=\frac{1}{\Dt(x)}\prod_{j=1}^n(1+vx_j(1-uT_{t,x_j}))\Dt(x), 
\tag{3.10}\\
L(u,v)&=\frac{1}{\Dt(x)}\prod_{j=1}^n(T_{t,x_j}+vx_j(1-uT_{t,x_j}))\Dt(x).
\endalign
$$
}
\noindent
If we take the coefficients of $v^n$, formula (3.9) recovers formula 
(3.1) for $D_x(u)$.  
The leading coefficient in $v$ of formula (3.10) also gives the well-known 
formula $D_x(u)=\Dt(x)^{-1}\prod_{i=1}^n(1-uT_{t,x_i})\,\Dt(x)$ 
when $q=t$. 
\par
Before the proof of Proposition 3.1, we will give a remark on 
the symbol of a $q$-difference operator. 
For a given $q$-difference operator 
$$
P(x;T_{q,x})=\sum_{\mu} a_{\mu}(x) T_{q,x}^\mu\in\K(x)[T_{q,x}^{\pm1}]
\tag{3.11}
$$
with rational coefficients, 
we define the {\it symbol} $p(x;\xi)$ of $P(x;T_{q,x})$ to be the function 
$$
p(x;\xi)=\sum_{\mu} a_{\mu}(x) \xi^\mu\in \K(x)[\xi^{\pm}]
\tag{3.12}
$$
in the commuting variables $(x;\xi)$. 
Then the action of $P(x;T_{q,x})$ on a function $\varphi(x)$ is 
recovered by 
$$
P(x;T_{q,x})\varphi(x)=p(x;T_{q,y})\varphi(y)|_{y=x},
\tag{3.13}
$$
by using the duplicated variables.
\Proofof{Proposition 3.1}{
We will treat only the case of $K(u,v)$, since the same argument gives 
the corresponding results for $L(u,v)$. 
Note first that the symbol of the operator $K_m(u)$ is given by
$$
\align
\operatorname{Symb}(K_m(u))
&=\frac{1}{\Dt(x)}
\sum\Sb |J|=m\endSb x_J
\sum_{I\subset J}(-u)^{|I|} T_{t,x}^I(\Dt(x)) \xi_I
\tag{3.14}\\
&=\frac{1}{\Dt(x)}
\sum\Sb |J|=m\endSb x_J
\prod_{j\in J}(1-uT_{t,x_j}\xi_j) \Dt(x).
\endalign
$$
Hence we have
$$
\align
\operatorname{Symb}(K(u,v))
&=\frac{1}{\Dt(x)} 
\sum\Sb J\endSb v^{|J|}x_J
\prod_{j\in J}(1-uT_{t,x_j}\xi_j) \Dt(x)
\tag{3.15}\\
&=\frac{1}{\Dt(x)}
\prod_{j=1}^n(1+vx_j(1-uT_{t,x_j}\xi_j))(\Dt(x)). 
\endalign
$$
When $q=t$, this already implies
$$
K(u,v)\varphi(x)=
\frac{1}{\Dt(x)}
\prod_{j=1}^n(1+vx_j(1-uT_{t,x_j}))(\Dt(x)\varphi(x)),
\tag{3.16}
$$
since the operators $x_j(1-uT_{t,x_j}\xi_j)$ ($j=1,\ldots,n$)
commute with each other. 
This proves formula (3.10) for $K(u,v)$.  
Returning to the general setting, we rewrite formula (3.15) as follows:
$$
\align
\operatorname{Symb}(K(u,v)) 
&=\frac{1}{\Dt(x)}\sum_{\sigma\in \frak{S}_n}\ep(\sigma)
\prod_{j=1}^n(1+vx_j(1-uT_{t,x_j}\xi_j))
(x_1^{\dt_{\sigma(1)}}\cdots x_n^{\dt_{\sigma(n)}})
\tag{3.17}\\
&=\frac{1}{\Dt(x)}\sum_{\sigma\in \frak{S}_n}\ep(\sigma)
\prod_{j=1}^n x_j^{\dt_{\sigma(j)}}
(1+vx_j(1-ut^{\dt_{\sigma(j)}}\xi_j))\\
&=\frac{1}{\Dt(x)}\det(x_j^{\dt_i}
(1+vx_j(1-ut^{\dt_i}\xi_j)))_{1\le i,j\le n}.
\endalign
$$ 
Consider now the determinant of Proposition:
$$
\align
&\frac{1}{\Dt(x)}\det(x_j^{\dt_i}(1+vx_j(1-ut^{\dt_i}T_{q,x_j})))_{1\le i,
j\le n}
\tag{3.18}\\
&=\frac{1}{\Dt(x)}\sum_{\sigma\in \frak{S}_n}\ep(\sigma)
\prod_{j=1}^n x_j^{\dt_{\sigma(j)}}
(1+vx_j(1-ut^{\dt_{\sigma(j)}}T_{q,x_j}))
\endalign
$$
Note that the factors in the product of the last expression are already 
ordered normally and that they act on different variables. 
{}From this we see that the symbol of the operator of (3.18) is precisely 
given by the formula (3.17).
This implies formula (3.9) of Proposition. 
}
Proposition 2.5 is obtained by expanding the determinantal formulas for 
$K(u,v)$ and $L(u,v)$ with respect to $v$. 
In fact we have
$$
\align
K(u,v)&=\sum_{w\in W} \ep(w) w\left(
x_1^{\dt_1}\cdots x_n^{\dt_n}
\prod_{j=1}^n(1+v x_j(1-u t^{\dt_j}T_{q,x_j}))
\right)
\tag{3.19}\\
&=\sum_{m=0}^n v^m \sum_{w\in W} \ep(w) w\left(
x_1^{\dt_1}\cdots x_n^{\dt_n}
e_m(X_1(u),\ldots,X_n(u)),
\right)
\endalign
$$
where we set
$$
X_j(u)=x_j(1-u t^{\dt_j}T_{q,x_j})\quad(j=1,\ldots,n). 
\tag{3.20}
$$
Namely we have 
$$
K_m(u)=\sum_{w\in W} \ep(w) w\left(
x_1^{\dt_1}\cdots x_n^{\dt_n}
e_m(X_1(u),\ldots,X_n(u))\right). 
\tag{3.21}
$$
Formula (2.10) for $K^+_m$ is the special case when $u=t^{m-n+1}$. 
A similar computation for $L(u,v)$ gives formula (2.10) for 
$K^-_m=t^{-\binom{n-m}{2}}L_m(t^{m-n+1})$. 
\section{\S4: Proof of Theorem 2.2}
The main idea in our proof of Theorem 2.2 is to use the generating function 
$$
\Pi(x,y)=\prod\Sb i\in[1,n]\\ j\in[1,m]\endSb
\frac{(tx_iy_j;q)_\infty}{(x_iy_j;q)_\infty}
=\sum_{\ell(\ld)\le m} b_\ld P_\ld(x) P_\ld(y) 
\tag{4.1}
$$
for the Macdonald polynomials in $x=(x_1,\cdots,x_n)$ 
and in $y=(y_1,\cdots,y_m)$, where $0\le m\le n$.
\Lemma{4.1}{
Let $K_m : \K[x]^W\to \K[x]^W $ be an operator acting on symmetric 
functions. 
Then $K_m$ has the property that 
$$
K_mJ_\ld(x)=J_{\ld+(1^m)}(x) \quad\text{for any  partition} \ \ 
\ld \ \ \text{with}\ \ \ell(\ld)\le m,
\tag{4.2}
$$
if and only if the following equality is verified:
$$
K_{m,x} \, \Pi(x,y)=\frac{1}{y_1\cdots y_m}D_y(1) \, \Pi(x,y),
\tag{4.3}
$$
where the suffix $x$ in $K_{m,x}$ indicates that $K_m$ should act 
on the $x$ variables. 
}
\Proof{
In fact, we have 
$$
\align
y_1\cdots y_m \, K_{m,x} \, \Pi(x,y) 
&=\sum_{\ell(\ld)\le m} b_\ld K_m (P_\ld(x)) \, y_1\cdots y_m \, P_\ld(y)
\tag{4.4}\\
&=\sum_{\ell(\ld)\le m} b_\ld K_m (P_\ld(x)) \,P_{\ld+(1^m)}(y).
\endalign
$$
As to the action of $D_y(1)$, we have 
$$
\align
&D_y(1)  \, \Pi(x,y) =\sum_{\ell(\ld)\le m} b_\ld \, P_\ld(x) \, 
d^m_{\ld}(1)\,P_\ld(y) \tag{4.5}\\
&\quad=\sum_{\ell(\ld)\le m} 
d^m_{\ld+(1^m)}(1) \, b_{\ld+(1^m)} \, P_{\ld+(1^m)}(x) \, P_{\ld+(1^m)}(y),
\endalign
$$
since 
$d^m_{\ld}(1)=(1-t^{m-1}q^{\ld_1})\cdots(1-tq^{\ld_{m-1}})(1-q^{\ld_m})=0$ 
if $\ld_m=0$. 
Hence, equality (4.3) is equivalent to  
$$
b_\ld\, K_m P_\ld(x) = d^m_{\ld+(1^m)}(1) \, b_{\ld+(1^m)}
\, P_{\ld+(1^m)}(x)\quad\text{for any}\ \ 
\ld \ \ \text{with}\ \ \ell(\ld)\le m.
\tag{4.6}
$$
This is equivalent to property (4.2) since $J_\ld(x)=c_\ld P_\ld(x)$ and 
$$
\frac{b_\ld}{c_\ld}=\frac{b_{\ld+(1^m)}}{c_{\ld+(1^m)}} d^m_{\ld+(1^m)}(1) 
\tag{4.7}
$$
by (1.4), (1.6) and (1.9). }
\par\medpagebreak
We now prove the equality (4.3) for our operators $K_m=K^\pm_m$ 
defined in (2.1). 
We will mainly consider the case of $K_m=K^+_m$: 
$$
K^+_{m,x} \, \Pi(x,y)=\frac{1}{y_1\cdots y_m} D_y(1)\,\Pi(x,y),
\tag{4.8}
$$
since the same argument is valid for $K^-_m$ as well. 
Firstly, the two sides of (4.8) can be written as follows:
$$
\align
& K^+_{m,x} \,\Pi(x,y)
\tag{4.9}\\
&=\Pi(x,y)
\sum\Sb J\subset[1,n]\\|J|=m\endSb x_J\sum_{I\subset J}
(-t^{m-n+1})^{|I|}\frac{T_{t,x}^I(\Dt(x))}{\Dt(x)}
\prod\Sb i\in I\\ k\in[1,m]\endSb\frac{1-x_iy_k}{1-tx_iy_k}
\endalign
$$
and
$$
\align
&\frac{1}{y_1\cdots y_m}D_y(1)\Pi(x,y)
\tag{4.10}\\
&=\Pi(x,y)
\frac{1}{y_1\cdots y_m}
\sum_{K\subset[1,m]}
(-1)^{|K|}\frac{T_{t,y}^K(\Dt(y))}{\Dt(y)}
\prod\Sb i\in[1,n]\\ k\in K\endSb\frac{1-x_iy_k}{1-tx_iy_k}.
\endalign
$$
Hence equality (4.8) is equivalent to
$$
\align
&\sum\Sb J\subset[1,n]\\|J|=m\endSb x_J\sum_{I\subset J}
(-t^{m-n+1})^{|I|}t^{\binom{|I|}{2}}
\prod\Sb i\in I\\ j\notin I \endSb
\frac{tx_i-x_j}{x_i-x_j}
\prod\Sb i\in I\\ k\in[1,m]\endSb\frac{1-x_iy_k}{1-tx_iy_k}
\tag{4.11}\\
&\quad=\frac{1}{y_1\cdots y_m}
\sum_{K\subset[1,m]}
(-1)^{|K|}t^{\binom{|K|}{2}}
\prod\Sb k\in K\\ \ell\notin K \endSb
\frac{ty_k-y_\ell}{y_k-y_\ell}
\prod\Sb i\in[1,n]\\ k\in K\endSb\frac{1-x_iy_k}{1-tx_iy_k}. 
\endalign
$$
The key point now is that this equality does {\it not} depend on $q$. 
It means that, in order to prove (4.8), we have only to establish it 
for some special value of $q$. 
\par
{}From now on, we consider the case of $q=t$, so that the Macdonald 
polynomials 
$P_\ld(x;q,t)$ reduce to the Schur functions 
$$
s_\ld(x)=\frac{\det(x_j^{\ld_i+\dt_i})_{1\le i,j\le n}}{\Dt(x)},
\quad \ld=(\ld_1,\cdots,\ld_n). 
\tag{4.12}
$$
Setting $q=t$, we have to prove that $K^+_m J_\ld(x;t,t)=J_{\ld+(1^m)}
(x;t,t)$, i.e., 
$$
K^+_m \, s_{\ld}(x)= s_{\ld+(1^m)}\prod_{k=1}^m (1-t^{\ld_k+m-k+1})
\tag{4.13}
$$
for any partition $\ld$ with $\ell(\ld)\le m$. 
It would then imply (4.8) for the case $q=t$ by Lemma 4.1 , 
hence (4.8) for the general $(q,t)$. 
With the notation of Section 3, we will compute the action 
of $K(u,v)$ on Schur functions by means of the formula (3.10):
$$
K(u,v)=\frac{1}{\Dt(x)}\prod_{j=1}^n(1+v x_j (1-u T_{q,x_j})) \Dt(x). 
\tag{4.14}
$$
In fact, we have
$$
\align
K(u,v) s_\ld(x)&=\frac{1}{\Dt(x)}
\prod_{k=1}^n(1+vx_k(1-uT_{t,x_k}))\det(x_j^{\ld_i+\dt_i})
\tag{4.15} \\
&=\frac{1}{\Dt(x)}\det(x_j^{\ld_i+\dt_i}(1+vx_j(1-ut^{\ld_i+\dt_i})))\\
&=\frac{1}{\Dt(x)}\sum_{K\subset[1,n]}v^{|K|}
\det(x_j^{\ld_i+\dt_i+\theta_K(i)}) \prod_{k\in K}(1-ut^{\ld_k+\dt_k}),
\endalign
$$
where $\theta_K(i)=1$ if $i\in K$ and $\theta_K(i)=0$ otherwise. 
Extending the notation (4.11) to any $n$-tuple of integers, 
we can rewrite (4.15) as 
$$
K(u,v) s_\ld(x) =\sum_{K\subset[1,n]} v^{|K|} 
s_{\ld+(1^K)}(x)  \prod_{k\in K}(1-ut^{\ld_k+\dt_k}),
\tag{4.16}
$$
where $(1^K)=(\theta_K(i))_{1\le i\le n}$. 
Suppose now that $\ld$ is a partition with $\ell(\ld)\le m$.
Since
$$
\ld+\dt=(\ld_1+n-1,\ldots,\ld_m+n-m,n-m-1,\ldots,0),
\tag{4.17}
$$
it is easily seen that $s_{\ld+(1^K)}(x)=0$ unless $K\subset[1,m]$ or 
$K\cap[m+1,n]=[m+1,m+r]$ for some $1\le r\le n-m$. 
We specialize the value of $u$ to be $t^{-(n-m-1)}$
so that $\prod_{k\in K}(1-ut^{\ld_k+\dt_k})=0$ if 
$K\cap[m+1,n]\ne\emptyset$.
Then we see
$$
K(t^{m-n+1},v) s_\ld(x)=\sum_{K\subset[1,m]} 
s_{\ld+(1^K)}(x)  \prod_{k\in K}(1-t^{\ld_k+m-k+1})
\tag{4.18}
$$
Namely we have 
$$
K_\ell(t^{m-n+1}) s_\ld(x)=\sum\Sb K\subset[1,m] \\ |K|=\ell\endSb 
s_{\ld+(1^K)}(x)  \prod_{k\in K}(1-t^{\ld_k+m-k+1})
\tag{4.19}
$$
for any partition $\ld$ with $\ell(\ld)\le m$. 
In particular we obtain
$$
\align
K_m(t^{m-n+1}) s_\ld(x)&=s_{\ld+(1^m)}\prod_{k=1}^m(1-t^{\ld_k+m-k+1}),
\tag{4.20}\\
K_\ell(t^{m-n+1}) s_\ld(x)&=0\quad\text{if}\quad\ell>m.
\endalign
$$
\par
This shows that, when $q=t$, $K^+_m=K_m(t^{m-n+1})$ has property (4.2).  
It implies (4.8) for the case $q=t$ by Lemma 4.1, hence 
(4.8) for the general $(q,t)$ via equality (4.11).  
Again by Lemma 4.1, we see that, for any $(q,t)$, 
the operator $K^+_m$ has the property (4.2) of raising operators for the 
Macdonald polynomials $J_\ld(x)$. 
This completes the proof of Theorem 2.2 for $K_m=K^+_m$. 
\par
The same argument is valid for $K_m=K^-_m$ as well. 
For comparison, we will include some formulas for the $q$-difference 
operators $K^-_m=t^{-\binom{n-m}{2}}L_m(t^{m-n+1})$ and 
$L(u,v)$. 
Formula (4.3) for $K^-_m$ is equivalent to
$$
\align
&\sum\Sb J\subset[1,n]\\|J|=m\endSb x_J\sum_{I\subset J}
(-t)^{m-|I|}t^{\binom{m-|I|}{2}}
\prod\Sb i\in I\\ j\notin I \endSb
\frac{x_i-tx_j}{x_i-x_j}
\prod\Sb j\notin I\\ k\in[1,m]\endSb\frac{1-x_jy_k}{1-tx_jy_k}
\tag{4.21}\\
&\quad=\frac{1}{y_1\cdots y_m}
\sum_{K\subset[1,m]}
(-1)^{|K|}t^{\binom{|K|}{2}}
\prod\Sb k\in K\\ \ell\notin K \endSb
\frac{ty_k-y_\ell}{y_k-y_\ell}
\prod\Sb i\in[1,n]\\ k\in K\endSb\frac{1-x_iy_k}{1-tx_iy_k}. 
\endalign
$$
It can be proved by computing the action of $L(u,v)$ with $q=t$ 
on Schur functions:
$$
L(u,v) s_\ld(x)=\sum_{K\subset[1,n]} s_{\ld+(1^K)}(x)
\prod_{k\in K}(1-ut^{\ld_k+\dt_k})\prod_{\ell\notin K} t^{\ld_\ell+\dt_\ell}.
\tag{4.22}
$$
\par\medpagebreak
As byproducts of this proof we obtain several interesting formulas. 
Consider the case $q=t$, and specialize (4.16) and (4.22) to $\ld=0$. 
Then one can easily see that these formulas reduce to 
$$
\align
K(u,v)1&=\sum_{m=0}^n v^m e_m(x) (ut^{n-m};t)_m,
\tag{4.23}\\
L(u,v)1&=\sum_{m=0}^n v^m e_m(x) (ut^{n-m};t)_m t^{\binom{n-m}{2}}.
\endalign
$$
\Proposition{4.2}{
For each $m=0,1,\cdots,m$, one has
$$
\align
\sum\Sb J\subset[1,n]\\|J|=m\endSb x_J\sum_{I\subset J}
(-u)^{|I|}t^{\binom{|I|}{2}}
\prod\Sb i\in I\\ j\notin I \endSb
\frac{tx_i-x_j}{x_i-x_j}
&=e_{m}(x)(ut^{n-m};t)_m,
\tag{4.24}\\
\sum\Sb J\subset[1,n]\\|J|=m\endSb x_J\sum_{I\subset J}
(-u)^{m-|I|}t^{\binom{n-|I|}{2}}
\prod\Sb i\in I\\ j\notin I \endSb
\frac{x_i-tx_j}{x_i-x_j}
&=e_{m}(x)(ut^{n-m};t)_m t^{\binom{n-m}{2}} .
\endalign
$$
In particular one has
$$
\align
\sum\Sb J\subset[1,n]\\|J|=m\endSb x_J\sum\Sb I\subset J\\ |I|=r\endSb
\prod\Sb i\in I\\ j\notin I \endSb
\frac{tx_i-x_j}{x_i-x_j}
&=t^{(n-m)r}\left[\matrix m\\r\endmatrix\right]_t e_{m}(x),
\tag{4.25}\\
\sum\Sb J\subset[1,n]\\|J|=m\endSb x_J\sum\Sb I\subset J\\ |I|=r \endSb
\prod\Sb i\in I\\ j\notin I \endSb
\frac{x_i-tx_j}{x_i-x_j}
&=\left[\matrix m\\r\endmatrix\right]_t e_{m}(x),
\endalign
$$
for each $r=0,1,\cdots,m$.
}
\noindent
Formulas (4.25) are obtained from (4.24) by taking the coefficients 
of $u^r$ or $u^{m-r}$. 
\par
Our raising operators $K^\pm_m$ for Macdonald polynomials $J_\ld(x;q,t)$ 
can be 
understood as a systematic generalization of formulas of this kind. 
\section{\S5: $q$-Difference lowering operators}
By using similar ideas, we can also construct $q$-difference lowering 
operators for Macdonald polynomials. 
\par\medpagebreak
For each $m=0,1,\cdots,n$, we define the $q$-difference operators 
$M^+_m$ and $M^-_m$ as follows:
$$
\align
M^+_m&=\sum\Sb J\subset[1,n]\\ |J|=m\endSb \prod_{j\in J}\frac{1}{x_j}
\sum_{I\subset J} (-1)^{|I|} t^{\binom{|I|}{2}}
\prod\Sb i\in I\\ j\in[1,n]\backslash I\endSb
\frac{tx_i-x_j}{x_i-x_j}
\prod_{i\in I}T_{q,x_i}, 
\tag{5.1}\\
M^-_m&=\sum\Sb J\subset[1,n]\\ |J|=m\endSb \prod_{j\in J}\frac{1}{x_j}
\sum_{I\subset J} (-t^{n-m})^{m-|I|} t^{\binom{m-|I|}{2}}
\prod\Sb i\in I\\ j\in[1,n]\backslash I\endSb
\frac{x_i-tx_j}{x_i-x_j}
\prod_{j\in[1,n]\backslash I}T_{q,x_j}. 
\endalign
$$
\par
\Theorem{5.1}{For each $m=0,1,\cdots,n$, the $q$-difference operators 
$M^+_m$ and $M^-_m$ preserve the ring $\K[x]^W$ of symmetric polynomials.
Furthermore they are lowering operators for Macdonald polynomials 
in the sense that 
$$
M^\pm_mJ_\ld(x)=a^m_\ld J_{\ld-(1^m)}(x), 
\quad a^m_\ld=\prod_{i=1}^m(1-t^{m-i}q^{\ld_i})(1-t^{n-i+1}q^{\ld_i-1}),
$$
for any partition $\ld$ with $\ell(\ld)\le m$.
In particular one has $M^\pm_m J_\ld(x)=0$ if $\ell(\ld)<m$. 
}
For the proof of Theorem 5.1, we prove the equality 
$$
M_{m,x} \Pi(x,y)=y_1\cdots y_m D_y(t^{n-m+1})\Pi(x,y),\quad M_m=M^\pm_m, 
\tag{5.2}
$$
for the auxiliary variables $y=(y_1,\ldots,y_m)$.  
By the same argument as in Lemma 4.1, one can prove that equality (5.2) 
implies Theorem 5.1.
Also, equality (5.2) is reduced to the special case when $q=t$. 
\par
We now introduce the following $q$-difference operators $M(u,v)$ and $N(u,v)$:
$$
M(u,v)=\sum_{m=0}^n v^m M_m(u),\quad N(u,v)=\sum_{m=0}^n v^m N_m(u),
\tag{5.3}
$$
where 
$$
\align
M_m(u)&=\sum\Sb J\subset[1,n]\\ |J|=m\endSb \frac{1}{x_J}
\sum_{I\subset J} (-u)^{|I|} 
\frac{T_{t,x}^I(\Dt(x))}{\Dt(x)}
T_{q,x}^I, 
\tag{5.4}\\
N_m(u)&=\sum\Sb J\subset[1,n]\\ |J|=m\endSb \frac{1}{x_J}
\sum_{I\subset J} (-u)^{m-|I|} 
\frac{T_{t,x}^{I^{\text{C}}}(\Dt(x))}{\Dt(x)} 
T_{q,x}^{I^{\text{C}}},
\endalign
$$
so that $M^+_m=M_m(1)$, and $M^-_m=t^{-\binom{n-m}{2}} N_m(1)$. 
These operators have the determinantal formulas
$$
\align
M(u,v)&=\frac{1}{\Dt(x)}
\det(x_j^{\dt_i}(1+\frac{v}{x_j}(1-ut^{\dt_i}T_{q,x_j}))
)_{1\le i,j\le n},
\tag{5.5}\\
N(u,v)&=\frac{1}{\Dt(x)}
\det(x_j^{\dt_i}(t^{\dt_i}T_{q,x_j}+\frac{v}{x_j}(1-ut^{\dt_i}T_{q,x_j}))
)_{1\le i,j\le n}. 
\endalign
$$
\par
When $q=t$, these formulas reduce to 
$$
\align
M(u,v)&=\frac{1}{\Dt(x)}\prod_{j=1}^n(1+\frac{v}{x_j}(1-uT_{t,x_j})) \, 
\Dt(x), 
\tag{5.6}\\
N(u,v)&=\frac{1}{\Dt(x)}\prod_{j=1}^n(T_{t,x_j}+\frac{v}{x_j}
(1-uT_{t,x_j})) \, \Dt(x). 
\endalign
$$
Furthermore, their action on Schur functions are computed as follows:
$$
\align
M(u,v)s_{\ld}(x)&=\sum_{K\subset[1,n]} v^{|K|} 
s_{\ld-(1^K)}(x)\prod_{k\in K}(1-ut^{\ld_k+\dt_k}),
\tag{5.7}\\
N(u,v)s_{\ld}(x)&=\sum_{K\subset[1,n]} v^{|K|}
s_{\ld-(1^K)}(x)\prod_{k\in K}(1-ut^{\ld_k+\dt_k})\prod_{k\notin K}t^{\ld_k+
\dt_k}. 
\endalign
$$
If $\ell(\ld)\le m$ and $u=1$, we have
$$
\align
M_m(1)s_\ld(x)&=s_{\ld-(1^m)}(x) \prod_{k=1}^m(1-t^{\ld_k+n-k}),
\tag{5.8}\\
N_m(1)s_\ld(x)&=s_{\ld-(1^m)}(x) \prod_{k=1}^m(1-t^{\ld_k+n-k}) 
t^{\binom{n-m}{2}}.
\endalign
$$
This implies that $M^+_m=M_m(1)$ and $M^-_m=t^{-\binom{n-m}{2}}N_m(1)$
satisfy the equality (5.2) for the case $q=t$. 
As we already explained, it completes the proof of Theorem 5.1. 
\par\medpagebreak
\Proposition{5.2}{
For each $m=0,1,\ldots,n$, we have
$$
M^+_m=\frac{1}{\Dt(x)}\sum_{w\in W}\ep(w)w\left(
x_1^{\dt_1}\cdots x_n^{\dt_n}e_{m}
(\Xi_1,\ldots,\Xi_n)\right) 
\tag{5.9}
$$
where 
$$
\Xi_i=\frac{1}{x_i}(1-t^{n-i}T_{q,x_i}) \quad (i=1,\ldots,n).
\tag{5.10}
$$
Similarly we have
$$
M^-_m
=\frac{t^{\binom{n}{2}-\binom{n-m}{2}}}{\Dt(x)}\sum_{w\in W}\ep(w)w\left(
x_1^{\dt_1}\cdots x_n^{\dt_n}e_{m}
(\Xi'_1,\ldots,\Xi'_n)\right) T_{q,x_1}\cdots T_{q,x_n}
\tag{5.11}
$$
where $\Xi'_i= \Xi_it^{-\dt_i}T_{q,x_i}^{-1}$ $(i=1,\cdots,n)$.
}
{}From these expressions, we obtain the following quasi-classical limit. 
$$
{\Cal M}_m=\frac{1}{\Dt(x)}\sum_{w\in W}\ep(w)w\left(
x_1^{\dt_1}\cdots x_n^{\dt_n}e_{m}
(\Cal D_1,\ldots,\Cal D_n)\right) 
\tag{5.12}
$$
where 
$$
\Cal D_i=\frac{1}{x_i}( \alpha x_i\frac{\partial\ }{\partial x_i}+n-i) 
\quad (i=1,\ldots,n).
\tag{5.13}
$$
These differential operators are the lowering operators for Jack polynomials:
$$
\Cal M_m J^{(\alpha)}_\ld(x)= 
\prod_{i=1}^m (\alpha\ld_i+m-i)(\alpha(\ld_i-1)+n-i+1) 
J^{(\alpha)}_{\ld-(1^m)}(x),
\tag{5.14}
$$
for any partition $\ld$ with $\ell(\ld)\le m$. 
\par
\Remark{5.3}{
It would be an interesting problem to describe the algebra generated by the 
raising and lowering operators $K^\pm_m, M^\pm_m$ 
($m=0,1,\cdots,n$ or $\infty$). 
Another interesting problem is to construct the analogues of raising 
and lowering operators for Koornwinder's multivariable Askey-Wilson 
polynomials, as well as for their nonsymmetric versions. 
We are currently studying these problems and hope to report on them in 
the near future. 
}


\Refs
\widestnumber\key{DJK}
\ref\key{KN} \by A.N.\,Kirillov and M.\,Noumi
\paper Affine Hecke algebras and raising operators for Macdonald polynomials
\paperinfo preprint
\yr February 1996
\endref
%
\ref\key{Ma} \by I.G.\,Macdonald 
\book Symmetric Functions and Hall Polynomials {\rm (Second Edition)}
\bookinfo Oxford Mathematical Monographs
\publ Oxford University Press Inc., New York
\yr 1995
\endref
%
\ref\key{MN} \by K.\,Mimachi and M.\,Noumi
\paper Notes on eigenfunctions for Macdonald's 
$q$-difference operators
\paperinfo preprint \yr January,1996
\endref
%
\ref\key{LV1} \by L.\,Lapointe and L.\,Vinet
\paper A Rodrigues formula for the Jack polynomials 
and the Mac\-donald-Stanley conjecture
\paperinfo preprint CRM-2294
\yr 1995
\endref
%
\ref\key{LV2} \by L.\,Lapointe and L.\,Vinet
\paper Exact operator solution of the Calogero-Sutherland model
\paperinfo pre\-print CRM-2272
\yr 1995
\endref
%
\endRefs
\enddocument